\documentclass[12pt]{article}
\usepackage{amssymb}
\usepackage{amsmath}
\usepackage{amsfonts}

\begin{document}

\begin{center}
{\Large\sc Local cohomology, arrangements of subspaces and
monomial ideals.}\\

$\ $ \\ {by Josep \`Alvarez Montaner \footnote{Partially
supported by the University of Nice.}, Ricardo Garc\'{\i}a
L\'opez \footnote{Partially supported by the DGCYT PB98-1185
and by INTAS 97 1644.}\\  and Santiago Zarzuela Armengou
\footnote{Partially supported by the DGCYT PB97-0893}.}
\end{center}

\newcommand{\ea}{{\mathbb A}^{n}_{k}}
\newcommand{\dli}{\mbox{indlim}_{\,P}\,}
\newcommand{\dlii}{\mbox{indlim}^{(i)}_{\,P}\,}
\newcommand{\ili}{projlim}
\newcommand{\supp}{\mbox{supp}_{+}\,}
\newcommand{\cpx}{{\mathbb C}}
\newcommand{\cp}{{\mathbb C}^n}
\newcommand{\NN}{{\mathbb N}}
\newcommand{\h}{\mbox{h}}
\newcommand{\ZZ}{{\mathbb Z}}
\newcommand{\RR}{\mathbb R}
\newcommand{\Roos}{\mbox{Roos}}
\newcommand{\height}{\mbox{h}}
\newcommand{\Hom}{\mbox{Hom}}
\newcommand{\mc}{\mathcal}
\newcommand{\p}{\mathfrak{p}}
\newcommand{\q}{\mathfrak{q}}
\newcommand{\rr}{\mathfrak{r}}
\newcommand{\poplus}{\textstyle{\bigoplus}}
\newcommand{\mm}{\mathfrak{m}}
\newcommand{\QQ}{{\mathbb Q}}
\newcommand{\Ext}{\mbox{Ext}^1_R\,}
\newcommand{\Ker}{\mbox{Ker}\,}
\newcommand{\E}[1]{\mbox{E}_{R}(R/#1)}
\newcommand{\lc}[1]{H^{#1}_{I}(R)}
\newcommand{\R}{k[x_1,\ldots,x_{n}]}
\newcommand{\chc}[2]{H^{#1}_c(#2)}
\newcommand{\vs}{\vspace{6mm}}
\renewcommand{\labelenumi}{\roman{enumi})}

\vs \noindent {\bf 0.\ Introduction.} \vs

Let $\ea$ denote the affine space of dimension $n$ over a
field $k$, let $X\subset \ea$ be an arrangement of linear
subvarieties. Set $R=\R$ and let $I\subset R$ denote an
ideal which defines $X$. In this paper we study the local
cohomology modules
$
\lc{i} := indlim_{j}\mbox{Ext}_R^{i}\, (R/I^j,R),
$
with special regard of the case where the ideal $I$ is
generated by monomials. \vs

If $k$ is the field of complex numbers (or, more generally,
a field of characteristic zero), the module $\lc{i}$ is
known to have a module structure over the Weyl algebra
$A_n(k)$, and one can therefore consider its characteristic
cycle, denoted $CC(\lc{i})$ in this paper (see e.g.
\cite[I.1.8.5]{Bjo}). On the other hand, the arrangement $X$
defines a partially ordered set $P(X)$ whose elements
correspond to the intersections of irreducible components of
$X$ and where the order is given by inclusion. \vs

Our first result is the determination of the characteristic
cycles $CC(\lc{i})$ in terms of the cohomology of some
simplicial complexes attached to the poset $P(X)$. It
follows from the formulas obtained that, in either the
complex or the real case, these characteristic cycles
determine the Betti numbers of the complement of the
arrangement in $\ea$. In fact, it was proved by Goresky and
MacPherson that the Betti numbers of the complement of an
arrangement $X$ can be computed as a sum of non-negative
integers, one for each non-empty intersection of irreducible
components of $X$. These integers are dimensions of certain
Morse groups (cf. \cite[Part III, Theorems 1.3 and
3.5]{GM}). We will see that, over a field of characteristic
zero, one can give a purely algebraic interpretation of them
in terms of local cohomology. \vs

More precisely, in section 1, and following closely
Bj\"orner--Ekedahl's proof of the $\ell$-adic version of
Goresky--MacPherson' s formula, we will establish the
existence of a Mayer--Vietoris spectral sequence
\[
E_2^{-i,j}= \mbox{indlim}^{(i)}_{\,P(X)}\, H^j_{I_p}(R)
\Rightarrow H^{j-i}_{I}(R)
\]
where $p$ runs over $P(X)$, $I_p$ is the (radical) ideal of
definition of the irreducible variety corresponding to $p$,
and $\mbox{indlim}^{(i)}_{\,P(X)}\,$ is the $i$-th left
derived functor of the inductive limit functor in the
category of inductive systems of $R$--modules indexed by
$P(X)$ . The main ingredient in the proof is the
Matlis-Gabriel theorem on structure of injective modules.
Our formula for the characteristic cycle of local cohomology
follows essentially from the fact, proved in section 1 as
well, that this spectral sequence degenerates at the
$E_2$-term. \vs

In case the arrangement $X$ is defined by a monomial ideal
$I$, the local cohomology modules $H^i_I(R)$ have a natural
$\ZZ^n$-grading. In section 2 we relate the multiplicities
of the characteristic cycle of $H^i_I(R)$ to its graded
structure (proposition 2.1). Using results of M.
Musta\c{t}\u{a}, we can conclude that the the multiplicities
of the characteristic cycle of $H^i_I(R)$ determine and are
determined by the graded Betti numbers of the Alexander dual
ideal of $I$. The relation between the $\ZZ ^n$-graded
structure of $H^i_I(R)$ and the multiplicities of its
characteristic cycle has also been established by K.
Yanagawa in a recent pre-print (\cite{Ya3}).\vs

The degeneration of the Mayer--Vietoris spectral sequence
provides a filtration of each local cohomology module
$\lc{i}$, where the succesive quotients are given by the
$E_2$-term. In general, not all the extension problems
attached to this filtration have a trivial solution. This is
a major difference between the case we consider here and the
cases considered by Bj\"orner and Ekedahl. Namely, in the
analogous situation for the $\ell$-adic cohomology of an
arrangement defined over a finite field the extensions
appearing are trivial not only as extensions of $\mathbb Q
_{\ell}$-vector spaces but also as Galois representations,
and for the singular cohomology of a complex arrangement the
extensions appearing are trivial as extensions of mixed
Hodge structures (cf. \cite[pg. 179]{BE}). This contrasts
with the fact that the proof of the degeneration of the
Mayer-Vietoris spectral sequence for $\ell$-adic or singular
cohomology uses deeper facts than the proof of the
degeneration for local cohomology (it relies on the
strictness of Deligne's weight filtration). \vs

In section 3 we solve these extensions problems in case the
ideal $I$ is monomial (the result stated is actually more
general, in that we work in the category of
$\varepsilon$-straight modules, which is a slight variation
of a category introduced by K. Yanagawa in \cite{Ya2}, and
includes as objects the local cohomology modules considered
above). It turns out that these extensions can be described
by a finite set of linear maps, which for local cohomology
modules supported at monomial ideals can be effectively
computed in combinatorial terms from certain Stanley-Reisner
simplicial complexes (using the results in \cite{Mu}). \vs

If the base field is the field of complex numbers, then by
the Riemann-Hilbert correspondence the category of
$\varepsilon$-straight modules is equivalent to a full
subcategory of the category of perverse sheaves in $\cpx^n$
(with respect to the stratification given by the coordinate
hyperplanes). This category has been described in terms of
linear algebra by Galligo, Granger and Maisonobe in
\cite{GGM1}. Given such a perverse sheaf and a coordinate
hyperplane $H_i$, one can define its partial variation map
$v_i$ around $H_i$. In section 4, and using the description
in [loc. cit.], we prove that the category of
$\varepsilon$-straight modules is equivalent to the category
of those perverse sheaves such that $v_i=0$ for all $1\le i
\le n$ (equivalently, all partial monodromies are the
identity map). This allows to describe this last category
purely in terms of $\ZZ ^n$-graduations. \vs

If $(P, \le)$ is a poset, we will denote by $K(P)$ the
simplicial complex which has as vertexs the elements of $P$
and where a set of vertexs $p_0, \dots , p_r$ determines a
$r$-dimensional simplex if $p_0< \dots < p_r$. If $K$ is a
simplicial complex and $E$ is a $k$-vector space, we will
denote by $Simp_{\ast}(K;E)$ the complex of simplicial
chains of $K$ with coefficients in $E$. \vs

On dealing with arrangements defined over fields of
characteristic zero we will use some notions from $\cal
D$-module theory, we refer to \cite{Bo} or \cite{Bjo} for
unexplained terminology. All modules over a non-commutative
ring (or over a sheaf of non-commutative rings) will be
assumed to be left modules. On dealing with arrangements
defined over fields of positive characteristic we will refer
to the notion of $F$-module introduced in \cite[Definition
1.1]{LyF}. If $A$ is a ring, we denote by $\mbox{Mod}_A$ the
category of $A$-modules. We denote $\mbox{Vect}_k$ the
category of $k$-vector spaces. \vs

For arrangements of subspaces over a field of characteristic
zero defined by monomial ideals, an algorithm to compute the
characteristic cycles of the local cohomology modules
considered in this paper was given by the first named author
in \cite{Al}. \vs
\vs

\noindent {\bf 1.\ Filtrations on local cohomology modules.
} \vs

Let $S$ be an inductive system of $R$--modules. J.E. Roos
introduced in \cite{Roos} a complex which has as $i$-th
cohomology the $i$-th left derived functor of the inductive
limit functor evaluated at $S$ (and the dual notion for
projective systems as well, this is actually the case
treated by Roos in more detail). We recall his definition in
the case of interest for us: \vs

Let $(P,\le)$ be a partially ordered set, let $\cal C$ be an
abelian  category with enough projectives and such that the
direct sum functors are exact (usually, $\cal C$ will be a
category of modules, sometimes with enhanced structure: a
$\mathcal D$-module structure, a F-module structure or a
grading). We will regard $P$ as a small category which has
as objects the elements of $P$ and, given $p,q \in P$, there
is one morphism $p\to q$ if $p\le q$. A diagram over $P$ of
objects of the category $\cal C$ is by definition a
covariant functor $F:P \to \mathcal C$. Notice that the
image of $F$ is an inductive system of objects of $\mathcal
C$ indexed by $P$. The category which has as objects the
diagrams of objects of $\cal C$ and as functors the natural
transformations is abelian and will be denoted
Diag$(P,\mathcal C)$. \vs

\noindent{\bf Definition:} The Roos complex of $F$ is the
homological complex of objects of $\cal C$ defined by
\[
\Roos _k(F) := \poplus_{p_0<\dots <p_k}\,F_{p_0\dots p_k}\ ,
\]
where $F_{p_0\dots p_k}=F(p_0)$ and, if $i>0$ and we denote
by $\pi_{p_0 \dots \widehat{p_i}\dots p_k}$ the projection
from ${\bigoplus}_{p_0< \dots< p_k} F_{p_0 \dots p_k}$ onto
$F_{p_0 \dots \widehat{p_i}\dots p_k}$, the differential on
$F_{p_0\dots p_k}$ is given by
\[
F(p_0\to p_1) + \sum_{i=1}^k (-1)^i \pi_{p_0 \dots
\widehat{p_i}\dots p_k}.
\]
This construction defines a functor $\Roos _{\ast}(\cdot):
\mbox{Diag}(P,\cal C) \to \mbox{C}(\cal C)$, where
$\mbox{C}(\cal C)$ denotes the category of chain complexes
of objects of $\cal C$. It is easy to see that this functor
is exact and commutes with direct sums. \vs

Let $X\subset \ea$ be an arrangement defined by an ideal
$I\subset R$. Given $p\in P(X)$, we will denote by $X_p$ the
linear affine variety in $\ea$ corresponding to $p$ and by
$I_p\subset R$ the radical ideal which defines $X_p$ in
$\ea$. Notice that the poset $P(X)$ is isomorphic to the
poset of ideals $\{I_p\}_p$, ordered by reverse inclusion.
We denote by $h(p)$ the $k$-codimension of $X_p$ in $\ea$
(that is, $h(p)$ equals the height of the ideal $I_p$). The
height of an ideal $J\subset R$ will be denoted by $\h(J)$,
if $J_1\subseteq J_2$ are ideals of $R$, we set
$\h(J_2/J_1):=\h(J_2)-\h(J_1)$.  \vs

Let $M$ a $R$-module, $i\ge 0$ an integer.
Then one can define a diagram of $R$--modules $H^i
_{[\ast]}(R)$ on the poset $P(X)$ by
\[
H^i _{[\ast]}(M) : \ \ p \mapsto H^i_{I_p}(M).
\]
This defines a functor $H^i_{[\ast]}(\cdot): \mbox{Mod}_R
\to \mbox{Diag}(P(X), \mbox{Mod}_R)$. \vs

\noindent {\bf Lemma:}{\it\ If $E$ is an injective
$R$--module, then the augmented Roos complex
\[
\Roos_{\ast} (H^0_{[\ast]}(E)) \to  H^0_{I}(E) \to 0
\]
is exact.} \vs

\noindent {\it Proof }: Since both $\Roos_{\ast}(\cdot)$ and
$H^0_{[\ast]}(\cdot)$ commute with direct sums, by the
Matlis--Gabriel theorem we can assume that there is a prime ideal
$\p \subset R$ such that $E=E_R(R/\p)$, the injective envelope of
$R/\p$ in the category of $R$-modules. Notice also that for any
ideal $J\subset R$, $H^0_{J}(E_R(R/\p))=E_R(R/\p)$ if $\p\subseteq
J$ and is zero otherwise. It will be enough to prove that if $\mm
\subset R$ is a maximal ideal, then the complex
\[
(\Roos_{\ast} (H^0_{[\ast]}(E)))_{\mm} \to
(H^0_{I}(E))_{\mm} \to 0
\]
is exact. If $I \not\subset \mm$, this complex is zero. Otherwise,
it equals the augmented complex
$Simp_{\ast}(K,E_{R_{\mm}}(R_{\mm}/\p R_{\mm}))\to
E_{R_{\mm}}(R_{\mm}/\p R_{\mm})\to 0$, where $K$ is the simplicial
complex attached to the subposet of $P(X)$ which has as vertexs
those linear subspaces $X_p$ such that $I_p\subseteq \mm$. It is
easy to see that $K$ is contractible, and then the lemma follows.
$\square$ \vs

Fix now an injective resolution $0\to R \to E^{\ast}$ of $R$
in the category of $R$--modules. Each of the modules $E^j\
(j\ge 0)$ defines a diagram $H^0_{[\ast]}(E^j)$ over $P(X)$
and one obtains a double complex
\[
\Roos _{-i}(H^0 _{[\ast ]} (E^j)), \ \ \ i\le 0\ , \ j\ge 0.
\]
(the change of sign on the indexing of the Roos complex is
because we prefer to work with a double complex which is
cohomological in both degrees). This is a second quadrant
double complex with only a finite number of non-zero
columns, so it gives rise to a spectral sequence that
converges to $H^{\ast}_{I}(R)$ (because of the lemma above).
More precisely, we have
\[
E_1 ^{-i,j} = \Roos_{i}(H^{0}_{[\ast]} (E^j)) \Rightarrow
H^{j-i}_{I}(R).
\]
The differential $d_1$ is that of the Roos complex, since
this complex computes the $i$-th left derived functor of the
inductive limit, the $E_2$ term will be
\begin{equation}
E_2^{-i,j} = \dlii H^j_{[\ast]}(R) \Rightarrow
H^{j-i}_{I}(R)
\end{equation}
(we write $P$ for $P(X)$ in order to simplify the writing of
our formulas). Hereafter this sequence will be called
Mayer-Vietoris spectral sequence. \vs

\noindent (1.1) {\it Remark:}\ If the base field $k$ is of
characteristic zero, we can choose an injective resolution
of $R$ in the category of modules over the Weyl algebra
$A_n(k)$. Since $A_n(k)$ is free as an $R$-module, it
follows (see e.g. \cite[II.2.1.2]{Bjo}) that this is also an
injective resolution of $R$ in the category of $R$-modules.
If the base field $k$ is of characteristic $p>0$, the ring
$R$ has a natural $F$-module structure and its minimal
injective resolution is a complex of $F$-modules and
$F$-module homomorphisms (see \cite[(1.2.b'')]{LyF}).
Therefore, the spectral sequence above may be regarded as a
spectral sequence in the category of $A_n(k)$--modules
(respectively, of $F$-modules). The main result of this
section is the following: \vs

\noindent (1.2) {\bf Theorem:} \ {\it Let $X\subset \ea$ be
an arrangement of linear varieties. Let $K(>p)$ be the
simplicial complex attached to the subposet $\{ q\in P(X)
\mid q > p\}$ of $P(X)$. Then:
\begin{enumerate}
\item There are $R$-module isomorphisms
$$\dlii H^j_{[\ast]}(R) \simeq \poplus_{h(p)=j}
\,(H^j_{I_p}(R))^{\oplus\, b_{i,p}}, $$ where
$b_{i,p}:=\dim_k \tilde{H}_{i-1}(K(>p); k)$, and $\tilde{H}$
denotes reduced simplicial homology. We agree that
$\tilde{H}_{-1}(K(>p); M)=M$, and that the reduced homology
of the empty simplicial complex is $M$ in degree $-1$ and
zero otherwise.
\item The Mayer--Vietoris spectral
sequence
\[
E_{2}^{-i,j} = \dlii H^j_{[\ast]}(R) \Rightarrow
H^{j-i}_{I}(R)
\]
degenerates at the $E_2$--term.
\item If $k$ is a field of characteristic zero, the isomorphisms in
i) are also isomorphisms of $A_n(k)$--modules and the spectral
sequence in ii) is a spectral sequence of $A_n(k)$--modules.
If $k$ is a field of positive characteristic, the
isomorphisms in i) are also isomorphisms of $F$--modules and
the spectral sequence in ii) is a spectral sequence of
$F$--modules.
\end{enumerate}
} \vs

\noindent {\it Proof of i)}\ (c.f. \cite[Proposition
4.5]{BE}): Given $p\in P(X)$ and a $R$--module $M$ we
consider the following three diagrams
\begin{description}
\item $F_{M,\ge p}$\ , defined by $F_{M,\ge p}(q)=M$ if $q\ge p$ and
$F_{M,\ge p}(q)=0$ otherwise.
\item $F_{M,> p}$\ , defined by $F_{M,> p}(q)=M$ if $q> p$ and
$F_{M,> p}(q)=0$ otherwise.
\item $F_{M,\,p}$\ , defined by $F_{M,\,p}(q)=M$ if $q=p$ and
$F_{M,\, p}(q)=0$ otherwise
\end{description}
(in all three cases $F(p\to q)=id$ if $F(p)=F(q)$ and it is
zero otherwise). In the category of diagrams of $R$--modules
over $P(X)$ we have an exact sequence
\[
0\to F_{M,\,>p} \to F_{M,\,\ge p} \to F_{M,\, p} \to 0
\]
Let $K(\ge p)$ be the simplicial complex attached to the
subposet $\{ q\in P(X) \mid q\ge p\}$ of $P(X)$. Then one
has
\[
\Roos_{\ast}(F_{M,\ge p}) = \mbox{Simp}_{\ast}(K(\ge p), M)
\ \ \mbox{ and } \ \ \Roos_{\ast}(F_{M,> p}) =
\mbox{Simp}_{\ast}(K(> p), M).
\]
Since the complex $K(\ge p)$ is contractible (to the vertex
corresponding to $p$), the long exact homology sequence
obtained from the sequence of complexes
\[
0\to \Roos_{\ast}(F_{M,\,>p}) \to \Roos_{\ast}(F_{M,\,\ge
p}) \to \Roos_{\ast}(F_{M,\, p}) \to 0
\]
gives
\[
\dlii F_{M,\,p} \cong \tilde{H}_{i-1}(K(>p); M),
\]
where the tilde denotes reduced homology, we agree that
$\tilde{H}_{-1}(K(>p); M)=M$, and
the reduced homology of the empty simplicial complex is $M$ in degree
$-1$ and zero otherwise. \vs

\noindent Notice that one has an isomorphism of diagrams
$H^j_{[\ast]}(R) \simeq
\poplus_{h(p)=j}\,F_{H^j_{I_p}(R),\,p}$. Thus,
\[
\dlii H^j_{[\ast]}(R) \cong \poplus_{h(p)=j}\,\dlii
F_{H^j_{I_p}(R),\,p} \cong \poplus_{h(p)=j}\,
\tilde{H}_{i-1}(K(>p); H^j_{I_p}(R))
\]
By the universal coefficient theorem, for $i>0$
\[
\tilde{H}_{i-1}(K(>p); H^j_{I_p}(R))\cong
(H^j_{I_p}(R))^{\oplus \, b_{i,p}}
\]
where $b_{i,p}:=\dim_k \tilde{H}_{i-1}(K(>p); k)$. Although
the isomorphism given by the universal coefficient theorem
is a priori only an isomorphism of $k$-vector spaces, it is
easy to check that in our case is also an isomorphism of
$A_n(k)$--modules (if char$(k)=0$) and of $F$-modules (if
char$(k)>0$). In particular, it is always an isomorphism of
$R$--modules. \vs

\noindent {\it Proof of ii) and iii):} Observe first that if
$I\subset R$ is an ideal and we set $h=\height(I)$, then all
associated primes of $H^h_{I}(R)$ are minimal primes of $I$
(this is well-known, we briefly sketch a proof: Let $\p$ be
an associated prime of $H^h_{I}(R)$. We have an spectral
sequence
\[
H^i_{\p R_{\p}}(H^j_{IR_{\p}}(R_{\p})) \Rightarrow
H^{i+j}_{\p R_{\p}}(R_{\p}).
\]
By assumption, $H^0_{\p R_{\p}}(H^h_{IR_{\p}}(R_{\p}))\neq
0$. Using that $H^{k}_{\p R_{\p}}(R_{\p})=0$ for $k\neq\dim
(R_{\p})$, one deduces that $\dim (R_{\p}/ IR_{\p}) = 0$).
\vs

It follows that if $\p,\q\subset R$ are prime ideals such that
$\p\not\subset \q$ and we set $i=\height(\p)$, $j=\height(\q)$,
then $\Hom_R(H^i_{\p}(R),H^j_{\q}(R))=0$. From this last fact and
i) above follows that the Mayer--Vietoris sequence degenerates at
the $E_2$--term. Part iii) follows from remark (1.1) above and the
observation at the end of the proof of part i). $\square$\vs

\noindent (1.3) {\bf Corollary:} \ {\it Let $X\subset \ea$
be an arrangement of linear varieties, let $I\subset R$ be
an ideal defining $X$. Then, for all $r\ge 0$ there is a
filtration $\{F^r_j\}_{r\le j\le n}$ of $H^r_{I}(R)$ by
$R$--submodules such that
\[
F^r_j/F^r_{j-1} \cong
\poplus_{h(p)=j}\,(H^j_{I_p}(R))^{\oplus\, m_{r,p}}.
\]
Once an injective resolution of $R$ as a $R$--module has been
fixed, this filtration is uniquely determined and functorial.
Moreover, if $char(k)=0$ the filtration can be chosen to be a
filtration by holonomic $A_n(k)$--modules and if char$(k)>0$ the
filtration can be chosen to be a filtration by $F$-modules. The
numbers $m_{r,p}$ can be computed in terms of the poset $P(X)$
attached to the arrangement as
\[
m_{r,p} = \dim _k \tilde{H}_{h(p)-r-1} (K(>p);k).
\]
If $k$ is a field of characteristic zero, the characteristic
cycle of the holonomic $A_n(k)$-module $H^r_{I}(R)$ is
\[
CC(H^r_I(R)) = \sum m_{r,p}\, T^{\ast} _{X_p} \ea\, ,
\]
where $T^{\ast} _{X_p} \ea$ denotes the relative conormal
subspace of $T^{\ast}\ea$ attached to $X_p$. \vs

\noindent If $k=\RR$ is the field of real numbers, the Betti
numbers of the complement of the arrangement $X$ in
${\mathbb A}^{n}_{\RR}$ can be computed in terms of the
multiplicities $\{ m_{i,p}\}$ as
\[
\dim_{\,\QQ} \tilde{H}_i({\mathbb A}^{n}_{\RR}-X;
\QQ)=\sum_{p}m_{i+1,\,p}.
\]
If $k=\cpx$ is the field of complex numbers, then one has
\[
\dim_{\,\QQ} \tilde{H}_i({\mathbb A}^{n}_{\cpx}-X;
\QQ)=\sum_{p}m_{i+1-h(p),\,p}.
\]} \vs

\noindent {\it Proof:} The filtration $\{F^r_j\}$ is the one
given by the degeneration of the Mayer--Vietoris spectral
sequence. It can be chosen to be a filtration by
$A_n(k)$-modules (respectively, $F$-modules) by (1.2.iii).
The formula for the characteristic cycle follows from the
fact that if $\height (I_p)=h$, then $CC(H^{h}_{I_p}(R)) =
T^{\ast}_{X_p}\,\ea$ and the additivity of the
characteristic cycle with respect to short exact sequences.
The formula for the Betti numbers of the complement
${\mathbb A}^{n}_{\RR}-X$ follows from a theorem of
Goresky--MacPherson (\cite[III.1.3. Theorem A]{GM}), which
states (slightly reformulated) that
\[
\tilde{H}_i({\mathbb A}^{n}_{\RR}-X; \ZZ) \cong \poplus_p
\,H^{h(p)-i-1}(K(\ge p),K(>p);\ZZ).
\]
Regarding a complex arrangement in ${\mathbb A}^{n}_{\mathbb C}$
as a real arrangement in ${\mathbb A}^{2n}_{\mathbb R}$, the
formula for the Betti numbers of the complement of a complex
arrangement follows from the formula for real arrangements.
$\square$\vs

\noindent (1.4) {\it Remarks:} i) Instead of an injective
resolution of $R$ in the category of $R$-modules, one could
take as well an acyclic resolution with respect to the
functors $\Gamma_J(\,\cdot\,)$, $J\subset R$ an ideal
(recall that $\Gamma _J(M) := \{ m\in M \mid \exists \,r \ge
0 \mbox{ such that } J^r m =0\}$). This fact will be used in
the next sections. \vs

ii)  Notice that the module $H^i_I(R)$ is only determined up
to unique isomorphism. The choice of an injective (or
acyclic) resolution of $R$ fixes an incarnation of
$H^i_I(R)$, and it is this incarnation which has been shown
to be filtered (however, if we choose another injective
resolution of $R$, then we have a filtered isomorphism
between the two incarnations of $H^i_I(R)$ associated to
both resolutions). \vs

iii) The formalism of Mayer--Vietoris sequences can be
applied to functors other than $H^{i}_{[\ast]}(\cdot)$ (and
other than those considered in \cite{BE}). For example, one
can consider the diagrams
\[
\mbox{Ext}^i_R\,(R/[\ast],R) : p \mapsto
\mbox{Ext}^i_R\,(R/I_p,R)
\]
and, similarly as for the local cohomology modules, one has
a spectral sequence
\[
E_2^{-i,j}=\dlii \mbox{Ext}^j_R\,(R/I_p,R) \Rightarrow
\mbox{Ext}^{j-i}_R\,(R/I,R)
\]
which degenerates at the $E_2$-term. Therefore, one can
endow the module $\mbox{Ext}^r_R\,(R/I,R)$ with a filtration
$\{G^r_j\}_{r\le j\le n}$ such that
\[
G^r_j/G^{r-1}_j \cong
\poplus_{h(p)=j}\,\mbox{Ext}^j_R\,(R/I_p,R)^{\oplus\,
m_{r,p}}.
\]
The functoriality of the construction gives that the natural
morphism \\ $\mbox{Ext}^r_R\,(R/I,R) \to H^r_I(R)$ is
filtered. In case $I\subset R$ is a monomial ideal, it would
be interesting to compare this filtration with the one
defined in \cite[Theorem 3.3]{Mu}. \newpage
\vs

 \noindent {\bf 2.\ Betti numbers vs. multiplicities.}
\vs

The ring $R= k[x_1, \dots , x_n]$ has a natural $\ZZ
^n$-graduation given by deg$(x_i)=\varepsilon _i$, where
$\varepsilon _1, \dots,\varepsilon _n$ denotes the canonical
basis of $\ZZ ^n$. If $M=\oplus_{\alpha \in \ZZ ^n}
M_{\alpha}$ and $N=\oplus_{\alpha \in \ZZ ^n} N_{\alpha}$
are graded $R$-modules, a morphism $f:M\rightarrow N$ is
said to be graded if $f(M_{\alpha}) \subseteq N_{\alpha}$
for all $\alpha \in \ZZ ^n$. Henceforth, the term graded
will always mean $\ZZ^n$-graded. We denote by
$^{\ast}\mbox{Mod}_{R}$ the category which has as objects
the graded $R$-modules and as morphisms the graded
morphisms. If $M,N$ are graded $R$-modules, we denote by
$^{\ast}\Hom _R(M,N)$ the group of graded morphisms from $M$
to $N$ (this group should not be confused with the internal
Hom in the category $^{\ast}\mbox{Mod}_R$, in particular it
is usually not a graded $R$-module). Its derived functors
will be denoted  $^{\ast}\mbox{Ext}_R^i(M,N),\ i\ge0$.\vs

We recall some facts about $^{\ast}\mbox{Mod}_{R}$ which will be
relevant for us (see \cite{GW} for proofs and related results): If
$M$ is a graded module one can define its $\ast$-injective
envelope $^{\ast}E(M)$ (in particular, $^{\ast}\mbox{Mod}_{R}$ is
a category with enough injectives). A graded version of the Matlis
-- Gabriel theorem holds: The indecomposable injective objects of
$^{\ast}\mbox{Mod}_{R}$ are the shifted injective envelopes
$^{\ast}E(R/\p)(\alpha)$, where $\p$ is an homogeneous prime ideal
of $R$ and $\alpha\in \ZZ ^n$, and every graded injective module
is isomorphic to a unique (up to order) direct sum of
indecomposable injectives. \vs

Injective objects of \,$^{\ast}\mbox{Mod}_{R}$ are usually
not injective as objects of \,$\mbox{Mod}_{R}$. However, if
$I\subseteq R$ is a monomial ideal, $\ast$-injective objects
are acyclic with respect to the functor
$\Gamma_I(\,\cdot\,)$. It follows that the local cohomology
modules $H^i_I(R)$ are objects of $^{\ast}\mbox{Mod}_R$. Our
aim is to relate the dimensions of its homogeneous
components to the multiplicities of its characteristic
cycle. \vs

A homogeneous prime ideal of $R$ is of the form $\p =
(x_{i_1},\dots,x_{i_k})$ \ $1\le i_1<\dots<i_k\le n$. We
will denote by $^{\ast}\mathcal{P}$ the set of homogeneous
prime ideals of $R$. Putting $\Omega =\{-1,0\}^n$, there is
a bijection
\begin{eqnarray*}
\Omega & \longrightarrow & ^{\ast}\mathcal{P}\\ \alpha &
\longmapsto & \ \p _{\alpha}\, =\ <x_{i_k} \, \mid \,
\alpha_{i_k} = -1 >.
\end{eqnarray*}

The preimage of an ideal $\p \in {^{\ast}\mathcal{P}}$ will
be denoted $\alpha_{\p}$ and, for $\alpha \in \Omega$, we
will set $\mid \alpha\mid = \mid \alpha_1\mid + \dots + \mid
\alpha_n \mid$. \vs

Let $I\subset R$ be a monomial ideal, $r\ge 0$ an integer,
let $X\subset \ea$ be the arrangement defined by $I$. It
will be convenient to reindex the multiplicities introduced
in section 1 as follows
\[
m _{r,\,\alpha} := \left\{
\begin{array}{l}
m_{r, p}\ \mbox{ if there is a } p\in P(X) \mbox{ with } \
\p_{\alpha}=I_p\,, \\ \  \\ 0 \ \ \ \ \mbox{ otherwise.}
\end{array} \right.
\]
Then we have \vs

\noindent (2.1) {\bf Proposition:} {\it In the situation and
with the notations described above
\[
m_{r,\,\alpha}= \,\dim _{\,k}H^r_I(R)_{\alpha} \ \ \mbox{
for all } \ r\ge 0\ ,\ \alpha\in \Omega.
\]}

\noindent{\it Proof:} Since $\ast$-injective modules are
$\Gamma_I$-acyclic, by (1.4,i) one can assume that the
filtration $\{F^r_j\}_{r\le j \le n}$ of $H^r_I(R)$
introduced in (1.3) has been constructed from an injective
resolution of $R$ in $^{\ast}\mbox{Mod}_R$. From this fact,
it follows that we can assume that the $R$-modules $F^r_j$
are graded, the inclusion maps $F^r_{j-1}\hookrightarrow
F^r_j$ are graded morphisms and one has isomorphisms of
graded modules
\[
F^r_j/F^r_{j-1}\cong \poplus_{\mid \alpha \mid
=j}\,(H^{j}_{\p_{\alpha}}(R))^{\oplus\, m_{r,\,\alpha}}.
\]
Notice that for $\beta \in \Omega$, one has
\[
(H^i_{\p _{\alpha}}(R))_{\beta}= \left\{
\begin{array}{l}
0 \ \ \mbox{if}\ \ \beta\neq\alpha\\ k \ \ \mbox{if} \ \
\beta = \alpha.
\end{array} \right.
\]
Using these facts and the exactness of the functors
\begin{minipage}{4.2cm}
\begin{eqnarray*}
^{\ast}\mbox{Mod}_R & \longrightarrow & \mbox{Vect}_k \ ,
\\ M & \longmapsto & \ M_{\alpha}
\end{eqnarray*}
\end{minipage}
the  desired result follows. $\square$\vs

Given a reduced monomial ideal $I\subset R$, the modules
$\mbox{Tor}^R_i(I,k)$ are $\ZZ^n$-graded and the graded
Betti numbers of $I$ are defined as
\[
\beta_{i,\,\alpha}(I):=\dim _k
\mbox{Tor}^R_i(I,k)_{\alpha}\, ,
\]
where $\alpha\in \ZZ ^n$. The Alexander dual ideal of $I$ is
the ideal
\[
I^{\vee}=\,<\ \Pi_{i\in S}\ x_i \mid S\subseteq
\{1,\dots,n\}\ ,\ \Pi_{i\not\in S}\ x_i\not\in I\ >\, .
\]
M. Musta\c{t}\u{a} showed in \cite{Mu} that if $\alpha\in
\Omega$, then one has
\[
\beta_{i,- \alpha}(I^{\vee})= \dim _k H^{\mid \alpha \mid -i
}_I(R)_{\alpha}\ ,
\]
and all other graded Betti numbers are zero. So, we have:
\vs

\noindent (2.2) {\bf Corollary:} {\it If $I$ is a reduced
monomial ideal and  $\alpha\in \{-1,0\}^n$, then
\[
\beta_{i,- \alpha}(I^{\vee})=m_{\mid \alpha \mid
-i,\,\alpha}.
\]}

In particular, if $\mbox{char}(k)=0$, the graded Betti
numbers of a monomial ideal $I$ can be obtained from the
($\mathcal D$-module theoretic) characteristic cycles of the
local cohomology of $R$ supported at its Alexander dual
$I^{\vee}$.\vs

Also, If $k = \cpx$ or $k=\RR$, it follows from the
corollary above and corollary (1.3) that the (topological)
Betti numbers of the complement in $\ea$ of the arrangement
defined by a monomial ideal $I$ can be obtained from the
(algebraic) graded Betti numbers of $I^{\vee}$. This fact
was already proved using a different approach in
\cite{GPW2}. \vs

\noindent {\bf 3.\ Extension problems.} \vs

If $M$ is a graded $R$-module and $\alpha\in\ZZ^n$, as usual
we denote by $M(\alpha)$ the graded $R$-module whose
underlying $R$-module structure is the same as that of $M$
and where the grading is given by $(M(\alpha))_{\beta}=
M_{\alpha + \beta}$. If $\alpha\in \ZZ^n$, we set
$x^{\alpha}=x_1^{\alpha _1}\cdots x_n^{\alpha _n}$ and $\supp
(\alpha)=\{ i \mid \,\alpha_i>0\,\}$. We recall the
following definition of K. Yanagawa: \vs

\noindent {\bf Definition {\rm (\cite[2.7]{Ya2})}:} A
$\ZZ^n$-graded module is said to be {\it straight} if the
following two conditions are satisfied:
\begin{enumerate}
\item $\dim_{\,k}M_{\alpha}<\infty $ for all $\alpha \in
\ZZ^n$.
\item The multiplication map $M_{\alpha}\ni y \mapsto
x^{\beta} y\in M_{\alpha + \beta}$ is bijective for all
$\alpha , \beta \in \ZZ^n$ with $\supp (\alpha + \beta ) =
\supp (\alpha)$.
\end{enumerate}
The full subcategory of the category $^{\ast}\mbox{Mod}_R$
which has as objects the straight modules will be denoted
{\bf Str}. Let $\varepsilon =- \sum_{i=1}^n \varepsilon_i =
(-1,\dots,-1)$. In order to avoid shiftings in local
cohomology modules, we will consider instead the following
(equivalent) category: \vs

\noindent {\bf Definition:} We will say that a graded module
$M$ is $\varepsilon$-straight if $M(\varepsilon)$ is
straight in the above sense. We denote $\varepsilon$-{\bf
Str} the full subcategory of $^{\ast}\mbox{Mod}_R$ which has
as objects the $\varepsilon$-{\bf Str} modules. It follows
from \cite{Mu}, \cite[Theorem 2.13]{Ya2} that if $I\subset
R$ is a monomial ideal and $r\ge 0$ is an integer,
$H^r_I(R)$ is an $\varepsilon$-straight module. \vs

\noindent (3.1) {\bf Proposition:\ }{\it Let $M$ be a
$\varepsilon$-straight module. There is a finite increasing
filtration $\{ F_j\}_{0\le j\le n}$ of $M$ by
$\varepsilon$-straight submodules such that for all $1\le
j\le n$ one has graded isomorphisms
\[
F_j/F_{j-1} \simeq \poplus_{\stackrel{\alpha \in
\Omega}{\scriptscriptstyle{ \mid \alpha \mid =
j}}}\,(H^j_{\p_{\alpha}}(R))^{\oplus\, m_{\alpha}}.
\]
where $m_{\alpha}=\dim_{\, k} M_{\alpha}$.} \vs

\noindent {\it Proof:} The existence of an increasing
filtration $\{G_j\}_j$ of $M$ by $\varepsilon$-straight
submodules such that all quotients $G_j/G_{j-1}$ are
isomorphic to local cohomology modules supported at
homogeneous prime ideals is an immediate transposition to
$\varepsilon$-straight modules of \cite[2.12]{Ya2}.
Inspection of Yanagawa's proof shows that, in order to prove
the existence of a filtration $\{ F_j\}_j$ satisfying the
condition of the proposition, it is enough to show that if
$\p _{\alpha}, \p_{\beta}$ are homogeneous prime ideals with
$\mid \alpha\mid = \mid \beta \mid =l$, then $^{\ast}\Ext
(H^l_{\p_{\alpha}}(R), H^l_{\p_{\beta}}(R))=0$. The minimal
$\ast$-injective resolution of $H^l_{\p_{\beta}}(R)$ is
\[
0 \to H^l_{\p_{\beta}}(R) \to
^{\ast}\E{\p_{\beta}}(-\varepsilon) \to \poplus_{i\mid \beta
_i=0} \,{^{\ast}\E{\p _{\beta -\varepsilon _i}}}
(-\varepsilon) \to \dots
\]
(see e.g. \cite[3.12]{Ya2}). Thus the vanishing of the above
$^{\ast}\mbox{Ext}$ module follows from the following: \vs

\noindent {\bf Claim:} For all homogeneous prime ideals
$\q\supset\p_{\beta}$ with $\mbox{h}(\q/\p_{\beta})=1$, one has
$^{\ast}\Hom_R(H^l_{\p_{\alpha}}(R),
{^{\ast}\E{\q}}(-\varepsilon))=0$.\vs

\noindent {\it Proof of the claim:} If $\q$ is such an homogeneous
prime ideal, we will assume that $\q\supset\p_{\alpha}$
(otherwise, the statement follows easily from the fact that $\q$
is the only associated prime of ${^{\ast}\E{\q}}$). Set $\alpha^c
:=\varepsilon - \alpha\in \Omega$. It suffices to prove the
vanishing of
$
^{\ast}\mbox{Hom}_R (H^l_{\p_{\alpha}}(R)(\varepsilon)\,,\,
{^{\ast}E}(R/\q))
$
By the equivalence of categories proved in \cite[2.8]{Ya2}, there
is a bijection between this group and
$
^{\ast}\mbox{Hom}_R(R/\p_{\alpha}(\alpha^c)\,,\,R/\q).
$
A graded morphism $\varphi: R/\p_{\alpha}(\alpha^c) \to
R/\q$ is determined by $\varphi(1)\in (R/\q)_{-\alpha^c}$.
But $(R/\q)_{-\alpha^c}=0$, so we are done. \vs

The equality $m_{\alpha}= \,\dim _{\,k}M_{\alpha}$ is proved
as in the proof of proposition (2.1). $\square$ \vs

\noindent (3.2) {\it Remark: } Even if $M$ is a local
cohomology module supported at a monomial ideal, in general
the submodules $F_j$ introduced in the proposition above are
not. This is one of the reasons for us to consider the
category of $\varepsilon$-straight modules.\vs

Henceforth we will assume that a $\varepsilon$-straight
module $M$ together with an increasing filtration $\{ F_j
\}_j$ of $M$ as in proposition (3.1) above have been fixed.
For each $j\le n$, one has then an exact sequence
\[
\hspace{-2cm}(s_j): \hspace{1cm} 0\to F_{j-1} \to F_{j}\to
F_j/F_{j-1} \to 0,
\]
which defines an element of $^{\ast}\Ext (F_j/F_{j-1}\,,
F_{j-1})$. In this section, our aim is to show that this
element is determined, in a sense that will be made precise
below, by the $k$-linear maps $\cdot x_i:M_{\alpha}\to
M_{\alpha+\varepsilon_i}$, where $\mid \alpha \mid =j$ and
$i$ is such that $\alpha_i=-1$. In particular, the sequence
$(s_j)$ splits if and only if $x_iM_{\alpha}=0$ for $\alpha,
i$ in this range. This is in general {\it not} the case for
local cohomology modules, which contrasts with the
situations considered in \cite{BE}, where the extensions
appearing are always split. We will prove first the
following lemma: \vs

\noindent {\bf Lemma: }{\it The natural maps}
\[
^{\ast}\Ext (F_j/F_{j-1}\,,F_{j-1}) \to
^{\ast}\Ext(F_j/F_{j-1}\,,F_{j-1}/F_{j-2})
\]
{\it are injective for all $j\ge 2$.} \vs

\noindent {\it Proof:} From the short exact sequence
$(s_{j-1})$, applying $^{\ast}\Hom (F_j/F_{j-1}, \, \cdot)$
we obtain the exact sequence
\[
^{\ast}\Ext(F_j/F_{j-1},F_{j-2}) \to ^{\ast}\Ext
(F_j/F_{j-1},F_{j-1}) \to
^{\ast}\Ext(F_j/F_{j-1},F_{j-1}/F_{j-2}),
\]
and we have to prove that
$^{\ast}\Ext(F_j/F_{j-1},F_{j-2})=0$. Applying again
$^{\ast}\Hom (F_j/F_{j-1},\, \cdot)$ to the exact sequences
$(s_l)$ for $l\le j-2$, and descending induction, the
assertion reduces to the statement $ ^{\ast}\Ext
(F_j/F_{j-1},F_l/F_{l-1})=0$ for $l\le j-2$. By proposition
$(3.1)$, it suffices to prove that if $\p, \q \subset R$ are
homogeneous prime ideals of heights $\height (\p) =j$ and
$\height(\q)=l$ and $l\le j-2$, then $^{\ast}\Ext
(H^j_{\p}(R),H^l_{\q}(R))=0$. As observed in the proof of
(3.1), the minimal $\ast$-injective resolution of
$H^l_{\q}(R)$ is
\[
0 \to H^l_{\q}(R) \to ^{\ast}\E{\q} (-\varepsilon)\to
\poplus_{\stackrel{h(\rr)=l+1}{\rr\in ^{\ast}\mathcal{P}}}
\,{^{\ast}\E{\rr}} (-\varepsilon)\to \dots
\]
Thus, again as in the proof of (3.1), it suffices to prove
that if $\height(\rr)=l+1<j$, then
$^{\ast}\Hom_R(H^j_{\p}(R), ^{\ast}\E{\rr}(\varepsilon))=0$.
This follows as in the proof of (1.2, ii), because $\rr$ is
the only associated prime of $\E{\rr}$ and $\p
\not\subset\rr$. $\square$ \vs

\noindent (3.3) {\it Remarks:} i) If $I\subset R$ is an ideal
defining an arbitrary arrangement of linear varieties in
$\ea$, one can prove an analogous lemma for the filtration
of the local cohomology module $H^r_I(R)$ introduced in
section 1.\vs

ii) The extension class of $(s_j)$ maps, via the morphism in
the above lemma, to the extension class of the sequence
\[
\hspace{-2cm}(s'_j): \hspace{1cm} 0\to F_{j-1}/F_{j-2} \to
F_{j}/F_{j-2}\to F_j/F_{j-1} \to 0.
\]

\noindent Let
\[
0\longrightarrow
F_{j-1}/F_{j-2}\longrightarrow{^{\ast}\mbox{E}}^0
\stackrel{d^0}{\longrightarrow}{^{\ast}\mbox{E}}^1
\longrightarrow\dots
\]
be the minimal $\ast$-injective resolution of
$F_{j-1}/F_{j-2}$. Given a graded morphism $F_j/F_{j-1}\to
\mbox{Im}\,d^0$, one obtains an extension of
$F_{j-1}/F_{j-2}$ by $F_j/F_{j-1}$ taking the following
pull-back: \vs

\begin{picture}(370,65)(-40,0)

\put(0,60){$0$}

\put(35,60){$F_{j-1}/F_{j-2}$}

\put(115,60){$^{\ast}\mbox{E}^0$}

\put(175,60){$\mbox{Im }d^0$}

\put(245,60){$0$}

\put(10,65){\vector(1,0){20}}

\put(85,65){\vector(1,0){25}}

\put(140,65){\vector(1,0){25}}

\put(215,65){\vector(1,0){25}}

\put(150,35){\dashbox{1}(9,9)}

\put(0,10){$0$}

\put(35,10){$F_{j-1}/F_{j-2}$}

\put(120,10){$E_{\varphi}$}

\put(175,10){$F_j/F_{j-1}$}

\put(185,35){$\varphi$}

\put(245,10){$0 \, ,$}

\put(10,15){\vector(1,0){20}}

\put(85,15){\vector(1,0){25}}

\put(140,15){\vector(1,0){25}}

\put(215,15){\vector(1,0){25}}

\put(60,35){$=$}

\put(56,25){\vector(0,1){27}}

\put(125,25){\vector(0,1){27}}

\put(180,25){\vector(0,1){27}}

\end{picture}
\vs

\noindent and all extensions of $F_{j-1}/F_{j-2}$ by
$F_j/F_{j-1}$ are obtained in this way. Take $\alpha \in
\Omega$ with $\mid \alpha \mid = j$. Applying the functor
$H^{\ast}_{\p_{\alpha}}(\,\cdot\,)$ to this diagram, we
obtain a commutative square \vs

\begin{picture}(165,80)(-80,0)
\put(0,0){$H^0_{\p_{\alpha}}(F_j/F_{j-1})$}

\put(0,70){$H^0_{\p_{\alpha}}(\mbox{Im}\, d^0)$}

\put(140,0){$H^1_{\p_{\alpha}}(F_{j-1}/F_{j-2})$}

\put(140,70){$H^1_{\p_{\alpha}}(F_{j-1}/F_{j-2})$}

\put(73,5){$\vector(1,0){62}$}

\put(60,75){$\vector(1,0){75}$}

\put(30,15){$\vector(0,1){45}$}

\put(160,15){$\vector(0,1){45}$}

\put(38,35){$\varphi^{\alpha}$}

\put(165,35){$\sim$}

\put(93,10){$\delta^{\alpha}$}

\put(90,80){$\sim$}
\end{picture}
\vs

\noindent where
\begin{enumerate}
\item The upper horizontal arrow is an isomorphism because
\\
${^{\ast}E^0} \simeq \oplus_{\mid \beta \mid =
j-1}{^{\ast}(\E{\p_{\beta}})}^{\oplus\, m_{\beta}}(-\varepsilon)$,
and then $H^i_{\p_{\alpha}}({^{\ast}E^0})=0$ for all $i\ge 0$.
\item The morphism $\delta^{\alpha}$ is the connecting
homomorphism of the given extension.
\end{enumerate}

Since $F_j/F_{j-1}= \oplus_{\mid \alpha\mid
=j}H^0_{\p_{\alpha}}(F_j/F_{j-1})$, the morphism $\varphi$
is determined by the morphisms $\varphi^{\alpha}$ for $\mid
\alpha \mid =j$, and these are in turn determined by the
connecting homomorphisms $\delta^{\alpha}$ via the
commutative square above. Observe also that the module
$H^1_{\p_{\alpha}}(F_j/F_{j-1})$ is isomorphic to  a direct
sum of local cohomology modules of the form
$H^j_{\p_{\alpha}}(R)$. We will show next that, because of
this fact, it suffices to consider the restriction of
$\delta^{\alpha}$ to the $k$-vector space of homogeneous
elements of multidegree $\alpha$.\vs

\noindent {\bf Lemma ({\rm linearization}):}\ {\it Let
$k_1,k_2\ge 0$ be integers, let
$M_1=H^l_{\p_{\alpha}}(R)^{\oplus\, k_1}$,
$M_2=H^l_{\p_{\alpha}}(R)^{\oplus\, k_2}$ . The restriction
map
\[
{^{\ast}\mbox{Hom}_R(M_1,M_2)}  \longrightarrow
\mbox{Hom}_{\,Vect_k}((M_1)_{\alpha},(M_2)_{\alpha})
\]
is a bijection. Moreover, these bijections are compatible
with composition of graded maps.} \vs

\noindent {\it Proof:}\ Taking the components of a graded
map $M_1 \to M_2$, it will enough to prove that all graded
endomorphisms of $N:=H^l_{\p_{\alpha}}(R)$ are
multiplications by constants of the base field $k$. Let
$\varphi: N\to N$ a graded endomorphism. Notice that
$\varphi$ can be regarded also as an endomorphism of the
graded module $N(\varepsilon)$. By \cite{Ya2},
$N(\varepsilon)$ is a straight module and $\varphi$ is
determined by its restriction to the $\NN ^n$-graded part of
$N(\varepsilon)$, which is $R/\p_{\alpha}(\alpha^c)$
($\alpha^c = \varepsilon - \alpha$). A graded $R$-module map
$R/\p_{\alpha}(\alpha^c)\to R/\p_{\alpha}(\alpha^c)$ is
determined by the image of $1 \in
(R/\p_{\alpha}(\alpha^c))_{-\alpha^c}$. Since
\[
(R/\p_{\alpha}(\alpha^c))_{-\alpha^c}=(R/\p_{\alpha})_0 = k,
\]
$\varphi$ must be the multiplication by some constant, as
was to be proved. $\square$ \vs

\noindent Thus, any extension of $F_{j-1}/F_{j-2}$ by
$F_j/F_{j-1}$ is determined by the $k$-linear maps
\[
\delta_{\alpha}^{\alpha}:
H^0_{\p_{\alpha}}(F_j/F_{j-1})_{\alpha} \longrightarrow
H^1_{\p_{\alpha}}(F_{j-1}/F_{j-2})_{\alpha}.
\]
One can easily check that
$H^0_{\p_{\alpha}}(F_j/F_{j-1})_{\alpha}\simeq M_{\alpha}$
and using a \u{C}ech complex one obtains that
\begin{eqnarray*}
\lefteqn{H^1_{\p_{\alpha}}(F_{j-1}/F_{j-2})_{\alpha} = \Ker
[\, \oplus_{\alpha_i = -1}
((F_{j-1}/F_{j-2})_{x_i})_{\alpha} \longrightarrow}
\\ & & \\ & &
\oplus_{\alpha_i = -1,\, \alpha_j=
-1}((F_{j-1}/F_{j-2})_{x_i x_j})_{\alpha}]
\stackrel{\sim}{\longrightarrow} \oplus_{\alpha_i=-1} \,
M_{\alpha + \varepsilon_i}
\end{eqnarray*}
where the last arrow is an isomorphism given by
multiplication by $x_i$ on
$((F_{j-1}/F_{j-2})_{x_i})_{\alpha}$. The connecting
homomorphism obtained applying $H^{\ast}_{\p_{\alpha}}(-)$
to the exact sequence $(s'_j)$ can be computed using
\u{C}ech complexes as well. It turns out that, via the
isomorphisms above, the corresponding map
$\delta_{\alpha}^{\alpha}$ is in this case precisely the map
\begin{eqnarray*}
M_{\alpha} & \longrightarrow & \oplus_{\alpha_i = -1}
M_{\alpha + \varepsilon_i} \\ m & \mapsto & \oplus (x_i
\cdot m).
\end{eqnarray*}
In conclusion, \vs

\noindent (3.4) {\bf Proposition:} {\it The  extension class
$(s_j)$ is uniquely determined by the $k$-linear maps $\cdot
x_i : M_{\alpha} \to M_{\alpha+\varepsilon_i} $ where $\mid
\alpha \mid = j$ and $\alpha_i = -1$.} \vs

\noindent (3.5) {\it Remark:} M. Musta\c{t}\u{a} has proved in
\cite{Mu} that for local cohomology modules supported at
monomial ideals, the linear maps $\cdot x_i :
H^j_I(R)_{\alpha}\to H^j_I(R)_{\alpha + \varepsilon_i}$ can
be explicitly computed in terms of the simplicial cohomology
of certain Stanley--Reisner complexes attached to $I$.\vs

\newpage

\noindent {\bf 4.\ Local cohomology and perverse
sheaves.} \vs

In this section we will use the following notations:
\begin{itemize}
\item $R = \cpx [x_1, \dots,x_n]$ (by a slight abuse of notation,
we will denote by $R$ as well the sheaf of regular algebraic
functions in $\cpx^n$).
\item $\mc O=$\, the sheaf of holomorphic
functions in $\cpx^n$.
\item $\mc D =$\, the sheaf of differential operators
in $\cpx^n$ with holomorphic coefficients.
\item $T=$ the union of the
coordinate hyperplanes in $\cpx^n$, endowed with the
stratification given by the intersections of its irreducible
components.
\item $\mc D^T_{hr}=$ the category of regular holonomic $\mc
D$-modules whose (singular) support is a union of
intersections of irreducible components of $T$.
\item $X_{\alpha}=$  the linear subvariety of $\cpx^n$
defined by the ideal $\p _{\alpha}\subset R$, \, $\alpha\in
\Omega$.
\end{itemize}\vs

\noindent By the Riemann-Hilbert correspondence, the functor
of solutions \\$\mathbb{R}\Hom_{\mathcal D}(-,\mathcal{O})$
establishes an equivalence between $\mathcal{D}_{hr}^T$, and
$Perv^T(\mathbb{C}^n)$, the category of complexes of sheaves
of finitely dimensional vector spaces on $\cpx^n$ which are
perverse relatively to the given stratification of $T$. In
\cite{GGM1}, the category $Perv^T(\mathbb{C}^n)$ has been
linearized as follows: Let $\mathcal{C}_n$ be the category
whose objects are families
$\{\mathcal{M}_{\alpha}\}_{\alpha\in \Omega}$ of finitely
dimensional complex vector spaces indexed by $\Omega := \{
-1, 0\}^n$, endowed with linear maps
\[
\mathcal{M}_{\alpha} \stackrel{u_i}{\longrightarrow}
\mathcal{M}_{\alpha - \varepsilon_i}\ \ , \ \
\mathcal{M}_{\alpha} \stackrel{v_i}{\longleftarrow}
\mathcal{M}_{\alpha - \varepsilon_i}\
\]
for each $\alpha\in\Omega$ such that $\alpha_i=0$. These
linear maps are required to satisfy the conditions
\[
u_iu_j = u_j u_i, \hskip 2mm v_i v_j =v_j v_i, \hskip 2mm
u_i v_j =v_j u_i \hskip 2mm \mbox{and} \hskip 2mm v_i u_i
+id \hskip 2mm \mbox{is invertible.}
\]
Such an object will be called an $n$-hypercube, the vector
spaces $\mathcal{M}_{\alpha}$ will be called its vertices
(for the description of the morphisms in this category, we
refer to \cite{GGM1}). It is proved in [loc. cit.] that
there is an equivalence of categories between
$Perv^T(\mathbb{C}^n)$ and $\mathcal{C}_n$. Given an object
$\mathcal{M}$ of $\mathcal{D}_{hr}^T$, the $n$-hypercube
corresponding to $\mathbb{R}\Hom_{\mathcal
D}(\mathcal{M},\mathcal{O})$ is constructed as follows: \vs

Consider $\mathbb{C}^n = \prod_{i=1}^n \mathbb{C}_i$, let
$K_i= \mathbb{R}^+ \subset \mathbb{C}_i$ and set
$V_i=\mathbb{C}_i \setminus K_i$. For any
$\alpha=(\alpha_1,\dots,\alpha_n) \in \Omega$ denote
\[
\mc{S}_{\alpha} := \frac{{\Gamma}_{\prod_{i=1}^n V_i}
\mathcal{O}} {\sum_{\alpha_k=-1} \Gamma_{\mathbb{C}_k \times
\prod_{i\neq k} V_i} \mathcal{O}}.
\]
Denoting with a subscript $0$ the stalk at the origin, one
has:
\begin{enumerate}
\item The vertices of the $n$-hypercube associated to $\mc M$
are the vector spaces
$\mathcal{M}_{\alpha}:=\Hom_{\mathcal{D}_0}(\mathcal{M}_0,
S_{\alpha,\, 0 })$.
\item The linear maps  $u_i$ are
those induced by the natural quotient maps
$\mathcal{S}_{\alpha} \to \mathcal{S}_{\alpha -
\varepsilon_i}$.
\item The linear maps $v_i$ are the
partial variation maps around the coordinate hyperplanes,
i.e. for any $\varphi \in
\Hom_{\mathcal{D}_0}(\mathcal{M}_0, \mathcal{S}_{\alpha\,, 0
})$ and $m\in \mathcal{M}_0$, $(v_i \varphi)(m) =
(\Phi_i(\varphi) - \varphi)(m)$, where $\Phi_i$ is the
monodromy around the hyperplane $x_i=0$.
\end{enumerate}
The following is proved as well in \cite{GGM1}, \cite{GGM2}:
\begin{enumerate}
\item[iv)] If $CC(\mathcal{M})=\sum
m_{\alpha} \hskip 1mm T_{X_\alpha}^{\ast} \mathbb{C}^n$ is
the characteristic cycle of $\mathcal{M}$, then for all
$\alpha\in\Omega$ one has the equality ${\dim}_{\mathbb{C}}
\mathcal{M}_{\alpha} = m_{\alpha}$.

\item[v)]
Let $\alpha, \beta \in\Omega$ be such that
$\alpha_i\beta_i=0$ for $1\le i\le n$. For each $j$ with
$\beta_j=-1$ choose any $\lambda_j \in \mathbb{C} \setminus
\mathbb{Z}$, set $\lambda_{\beta}=\{\lambda_j\}_j$ and let
$I_{\alpha,\beta,\lambda_{\beta}}$ denote the left ideal in
$\mathcal{D}$ generated by $ (\{x_i\,|\,\alpha_i=-1\},
\{\partial_k\,|\,\alpha_k = \beta_k = 0\}, \{x_j\partial_j -
\lambda_j\,|\,\beta_j=-1\})$. Then the simple objects of the
category $\mathcal{D}_{hr}^T$ are the quotients
$\frac{\mathcal{D}}{I_{\alpha,\beta,\lambda_{\beta}}}$.
\end{enumerate}

 \noindent {\bf Definition:}\, We say that an object $\mc
M$ of ${\mc D}_{hr}^T$ has variation zero if the morphisms
$v_i:\Hom_{{\mc D}_0}({\mc M}_0, \mc S_{\alpha -
{\varepsilon}_i,\, 0}) \longrightarrow \Hom_{{\mc D}_0}(\mc
M_0,\mc S_{\alpha,\,0})$ are zero for all $1\le i\le n$ and
all $\alpha\in\Omega$ with $\alpha_i=0$. It is easy to prove
that the modules with variation zero form a full abelian
subcategory of $\mathcal{D}_{hr}^T$, that will be denoted
${\mc D}^T_{0}$. \vs

\noindent (4.1) {\it Remarks:}\ a) One can check that if $f\in
\cpx[x_1, \dots, x_n]$ is a monomial, the $\mc D$-module
${\mc O}\,[1/f]$ has variation zero. Using a \v{C}ech
complex, it follows that  if $X$ is a subvariety of $\cpx^n$
defined by the vanishing of monomials, then the local
cohomology modules ${\mc H}^{\ast}_X(\mc O)$ are also
modules of variation zero. \vs

b) A simple object
$\frac{\mathcal{D}}{I_{\alpha,\beta,\lambda_{\beta}}}$ of
the category $\mathcal{D}_{hr}^T$ has variation zero if and
only if $\beta_k = 0$ for $1\le k\le n$. Thus, the simple
objects of $\mathcal{D}_{0}^T$ are of the form
\[
\frac{\mathcal{D}}{\mathcal{D}(\{x_i\,|\,\alpha_i=-1\},
\{\partial_j\,|\,\alpha_j = 0\})}.
\]
This module is isomorphic to the local cohomology module
${\mc H}_{X_{\alpha}}^{|\alpha|}(\mc O)$. Since every
holonomic module has finite length, we have: \vs

\noindent (4.2) \ {\bf Proposition:}\ {\it An object $\mc M$
of ${\mc D}^T_{hr}$ has variation zero if and only if there
exists a finite increasing filtration $0 = {\mc F}_0
\subseteq {\mc F}_1 \subseteq \cdots \subseteq {\mc F}_r =
{\mc M}$ by objects of ${\mc D}^T_{hr}$ such that the
quotients ${\mc F}_j/{\mc F}_{j-1}$ are isomorphic to local
cohomology modules ${\mc H}_{X_{\alpha}}^{|\alpha|}(\mc
O)$.} \vs

If $M$ is a $A_n(\cpx)$-module , then $\mc M ^{an}:= \mc O
\otimes_{R} M$ has a natural $\mc D$-module structure. This
allows to define a functor
\begin{eqnarray*}
(-)^{an}:\, Mod_{A_n(\cpx)} &\longrightarrow& Mod_{\mc D}.
\\ M &\longrightarrow & \mc M^{an} \\ f&\longrightarrow& id
\otimes f
\end{eqnarray*}
If $M$ is an $\varepsilon$-straight module, then $\mc
M^{an}$ is an object of $\mc D^T_{0}$. This follows from
(1.3), (4.1) and from the fact that, because the morphism $R
\to \mc O$ is flat, one has isomorphisms $\mc O \otimes _{R}
H ^l_{\p_{\alpha}}(R) \cong \mc H^l_{X_{\alpha}}(\mc O)$ for
all $\alpha\in\Omega$, $l\ge 0$. The main result of this
section is the following : \vs

\noindent (4.3) {\bf Theorem:}\,{\it The functor
\[
(-)^{an}: \varepsilon - \mathbf{ Str}\longrightarrow
\mathcal{D}_{0}^T
\]
is an equivalence of categories.} \vs

\noindent We will prove first the following lemma (which in
particular gives the fully faithfulness of $(-)^{an}$):\vs

\noindent (4.4) {\bf Lemma:} {\it Let $M, N$ be
$\varepsilon$-straight modules. For all \ $i\ge 0$, we have
functorial isomorphisms
\[
^{*}{\rm Ext}_{R}^i (M,N)\cong {\rm
Ext}_{\mathcal{D}^T_0}^i\, (\mc M^{an}\, ,\, \mc N^{an}).
\]}

\noindent {\it Proof of (4.4):}\, It has been proved by
Yanagawa that the category of straight modules has enough
injectives. It will follow from our proof that the category
$\mc D^T_0$ has enough injectives as well, so that both $\rm
Ext$ functors are defined and can be computed using
resolutions.\vs

By induction on the length we can suppose that $M$ and $N$
are simple objects, i.e. $M=H^r_{\p_{\alpha}}(R)$ and
$N=H^l_{\p_{\beta}}(R)$. Recall that the minimal
$\ast$-injective resolution of $N$ is:
\begin{equation}
0 \to H^l_{\p_{\beta}}(R) \to
^{\ast}\E{\p_{\beta}}(-\varepsilon) \to
\poplus_{i\mid\beta_i=0} \,{^{\ast}\E{\p_{\beta -
\varepsilon_i }}}(-\varepsilon) \to \cdots
\end{equation}
If $\alpha\in \Omega$, set $\mc E ^{\alpha}:=
(^{\ast}\E{\p_{\alpha}}(-\varepsilon))^{an}$. We claim that
for all $\alpha\in \Omega$, \, $\mc E^{\alpha}$ is an
injective object of $\mc D^T_0$. Using the description of
$\ast$-injective envelopes in \cite[3.1.5]{GW} (or,
alternatively, the equivalence of categories proved in
\cite[2.8]{Ya2}) one can see that, for all $\alpha \in
\Omega$, there are isomorphisms
\[
^{\ast}\E{\p_{\alpha}}(-\varepsilon) \cong \frac{R[
\frac{1}{x_1\cdots x_n}]}{\sum_{\alpha_i =-1} R[
\frac{1}{x_1 \cdots \widehat{x_i}\cdots x_n}]}.
\]
From this isomorphism it is easy to compute the
$n$-hypercube corresponding to $\mc E^{\alpha}$, namely one
has
\[
\mc E^{\alpha}_{\gamma}= \left\{
\begin{array}{l}
\cpx \ \mbox{ if } \gamma_i\le \alpha_i \mbox{ for all } \
1\le i\le n \,,
\\ \  \\ 0 \ \ \ \ \mbox{ otherwise.}
\end{array} \right.
\]

\noindent The map $u_i:\mc E^{\alpha}_{\gamma}\to \mc
E^{\alpha}_{\gamma - \varepsilon _i}$ is the identity if
$\mc E^{\alpha}_{\gamma }= \mc E^{\alpha}_{\gamma -
\varepsilon _i}=\cpx $ and it is zero otherwise. The
exactness of the functor $\mbox{Hom}_{\mc D}(-,\mc
E^{\alpha})$ on $\mc D^T_0$ can now be proved by passing to
the category $\mc C^n$ of $n$-hypercubes, where the
assertion reduces to a simple question of linear algebra
that is left to the reader. \vs

From flatness of $R \to \mc O$ and the injectivity of the
$\mc E ^{\alpha}$ proved above, it follows that one has the
following injective resolution of $\mc N^{an} = \mc
H^l_{X_{\beta}}(\mc O)$ in $\mc D^T_0$:
\begin{equation}
0\to \mc H^l_{X_{\beta}}(\mc O) \to \mc E^{\beta} \to
\poplus_{i\mid\beta_i=0}\, \mc E^{\beta - \varepsilon _i}
\to \cdots
\end{equation}
Let $K_1^{\bullet}$ be the complex obtained applying
$^{\ast}\mbox{Hom}_R(M,-)$ to the resolution (2) and let
$K_2^{\bullet}$ be the one obtained  applying
$\mbox{Hom}_{\mc D}(\mc H^r_{X_{\alpha}}(\mc O), -)$ to (3).
We have an injection $K_1^{\bullet}\hookrightarrow
K_2^{\bullet}$ and we want to show that it is an
isomorphism. We have
\begin{eqnarray}
^{*}{\rm Hom}_{R} (H^r_{\p_{\alpha}}(R),
^{\ast}\E{\p_{\gamma}}(-\varepsilon))= \left\{
\begin{array}{ll} \cpx & \mbox{ if } \alpha=\gamma
\\ 0 & \mbox{ otherwise }
\end{array} \right.
\end{eqnarray}
(this can be seen taking the positively graded parts and
using \cite[2.8]{Ya2}, as done before in similar
situations). The same equality holds replacing the left hand
side in (4) by ${\rm Hom}_{\mc D}(\mc H^r_{X_{\alpha}}(\mc
O)\, ,\, \mc E ^{\gamma})$. This is easily proved
considering the corresponding $n$-hypercubes, recall that
the $n$-hypercube corresponding to $\mc H^r_{X_{\alpha}}(\mc
O)$ is
\[
\mc H^r_{X_{\alpha}}(\mc O)_{\delta} = \left\{
\begin{array}{l}
\cpx \ \ \ \mbox{ if } \delta = \alpha
\\ \  \\ 0 \ \ \ \ \mbox{ otherwise.}
\end{array} \right.
\]
It follows that $K_1^{\bullet}\cong K_2^{\bullet}$, and then
we are done. \, $\square$ \vs

\noindent {\it Proof of (4.3):}\, By lemma (4.4) the functor
$(-)^{an}$ is fully faithful, so it remains to prove that it
is dense. Let $\mc N$ be an object of $\mc D^T_{0}$, let
$\mc N' \subseteq \mc N$ be a submodule such that $\mc
N'':=\mc N/\mc N'$ is simple. By induction on the length,
there are $\varepsilon$-straight $R$-modules $M'$ and $M''$
such that $\mc N'\cong (\mc M')^{an}$ and $\mc N''\cong (\mc
M'')^{an}$. The extension $0\to \mc N'\to \mc N\to \mc
N''\to 0$ correspond to an element $\xi$ of
$\mbox{Ext}^1_{\mc D ^T_0}(\mc N'',\mc N')$. Let $0\to M'\to
M\to M''\to 0$ be an extension such that its class in
$^{\ast}\mbox{Ext}^1_R(M'',M')$ maps to $\xi$ via the
isomorphism of lemma (4.4). One can check that $\mc N\cong
\mc M^{an}$, and then the theorem is proved. $\square$

\begin{footnotesize}

\end{footnotesize}

\bigskip

\begin{footnotesize}
Josep \`Alvarez Montaner

\smallskip
{\tt joalvarez@ma1.upc.es}

\smallskip
{\it Departament de Matem\`atica Aplicada I

Universitat Polit\`ecnica de Catalunya

Avinguda Diagonal, 647

Barcelona 08028, Spain }

\medskip

Ricardo Garc\'{\i}a L¢pez

\smallskip
{\tt rgarcia@mat.ub.es}

\medskip
Santiago Zarzuela Armengou

\smallskip
{\tt zarzuela@mat.ub.es}

\smallskip
{\it Departament d'\`Algebra i Geometria

Universitat de Barcelona

Gran Via, 585

Barcelona 08007, Spain }

\end{footnotesize}

\end{document}